\documentclass[11pt]{amsart}
\usepackage{amsmath}
\usepackage{amssymb}
\usepackage{bbm}
\usepackage{hyperref}
\usepackage[numbers,sort&compress]{natbib}
\usepackage{amscd}

\theoremstyle{definition}
\newtheorem{theorem}{Theorem}[section]
\newtheorem{lemma}[theorem]{Lemma}
\newtheorem{proposition}[theorem]{Proposition}
\newtheorem{definition}[theorem]{Definition}
\newtheorem{corollary}[theorem]{Corollary}
\newtheorem{criterion}[theorem]{Criterion}

\newtheorem{problem}[theorem]{Problem}

\begin{document}
\title{The classification problem for automorphisms of C*-algebras}
\author{Martino Lupini}
\address{Martino Lupini\\
Department of Mathematics and Statistics\\
N520 Ross, 4700 Keele Street\\
Toronto Ontario M3J 1P3, Canada, and Fields Institute for Research in
Mathematical Sciences\\
222 College Street\\
Toronto ON M5T 3J1, Canada.}
\email{mlupini@mathstat.yorku.ca}
\urladdr{http://www.lupini.org/}
\thanks{The author was supported by the York University Susan Mann
Dissertation Scholarship.}
\dedicatory{}

\begin{abstract}
We present an overview of the recent developments in the study of the classification problem for automorphisms of C*-algebras from the perspective of Borel complexity theory.
\end{abstract}

\maketitle
\tableofcontents

\section{Introduction\label{Section:Borel}}

Borel complexity theory is an area of logic that studies, generally
speaking, the relative complexity of classification problems in mathematics.
The main ideas and methods for such a study come from descriptive set
theory, which can be described as the analysis of definable sets and their
properties. In the framework of Borel complexity theory, a \emph{%
classification problem} is regarded as an equivalence relation on a Polish
space. Perhaps after a suitable parametrization, this virtually covers
almost all concrete classification problems in mathematics.

Under the assumption that the classes of objects under consideration are
naturally parametrized by the points of a Polish space, it is to be
expected---and demanded---that a satisfactory classification satisfy some 
\emph{constructibility assumption}. A sensible notion of constructibility is
to be \emph{Borel measurable} with respect to the given parametrizations.
This leads to the following definition, first introduced and studied in \cite%
{friedman_borel_1989}. Suppose that $E,E^{\prime }$ are equivalence
relations on standard Borel spaces $X,X^{\prime }$. A\ \emph{Borel reduction}
from $E$ to $E^{\prime }$ is a Borel function $f:X\rightarrow X^{\prime }$
such that%
\begin{equation*}
f\left( x_{1}\right) E^{\prime }f\left( x_{2}\right) \quad \text{if and only
if}\quad x_{1}Ex_{2}
\end{equation*}%
for every $x_{1},x_{2}\in X$. The function $f$ can be seen as a constructive
way to assign to the objects of $X$ complete invariants up to $E$ that are $%
E^{\prime }$-classes. Equivalently a Borel reduction can be seen as an
embedding of the quotient space $X/E$ into the quotient space $X^{\prime
}/E^{\prime }$ that is \textquotedblleft definable\textquotedblright , in
the sense that it admits a Borel lifting from $X$ to $X^{\prime }$.

The notion of Borel reduction\emph{\ }allows one to compare the complexity
of different classification problems. If there exists a Borel reduction from 
$E$ to $E^{\prime }$, then $E$ is \emph{Borel reducible }to $E^{\prime }$.
Such a notion formalizes the assertion that $E^{\prime }$ is \emph{at least
as complicated }as $E$ or, alternatively, that $E$ is at most as complicated
as $E^{\prime }$. In fact any Borel assignment of complete invariants up to $%
E^{\prime }$ yields a Borel assignment of complete invariants up to $E$ by
precomposing with a Borel reduction. Therefore if one knows how to classify
the objects of $X^{\prime } $ up to $E^{\prime }$, then one also knows via
the Borel reduction how to classify the objects of $X$ up to $E$.

Borel reducibility gives a way to compare the complexity of classification
problems. Some canonical equivalence relations are then used as \emph{%
benchmarks of complexity}. This allows one to build a hierarchy of
classification problems in mathematics. The first natural benchmark is
provided by the relation $=_{\mathbb{R}}$ of equality of real numbers. An
equivalence relation $E$ is \emph{smooth }if it is Borel reducible to $=_{%
\mathbb{R}}$. One can replace $\mathbb{R}$ with any uncountable standard
Borel space \cite[\S 15.6]{kechris_classical_1995}. Smoothness amounts at
saying that one can (explicitly) parametrize the equivalence classes of $E$
by the points of a standard Borel space. While the most satisfactory
classification results in mathematics are of this form, it was soon realized
that many natural and interesting equivalence relations are not smooth. This
observation goes back to the work of Mackey and Glimm on irreducible
representations of groups and C*-algebras \cite{mackey_groups_1957,
glimm_type_1961}.

One of the main goals of Borel complexity theory is to refine the analysis
in the nonsmooth case, and to understand if any satisfactory classification,
possibly involving a more general kind of invariants, is possible. A more
generous notion than smoothness is being \emph{classifiable by countable
structures}. In this case the benchmark of complexity is provided by the
relation $\cong _{\mathcal{L}}$ of isomorphism of $\mathcal{L}$-structures
for some first-order language $\mathcal{L}$. An equivalence relation is then
classifiable by countable structures if it is Borel reducible to $\cong _{%
\mathcal{L}}$ for some first-order language $\mathcal{L}$.

Assuming without loss of generality that $\mathcal{L}$ is a relational
language, the space of $\mathcal{L}$-structures can be regarded as the space 
$\prod_{R}2^{\mathbb{N}^{n_{R}}}$ where $R$ ranges over the symbols in $%
\mathcal{L}$ and $n_{R}$ denotes the arity of $R$. In this parametrization
the relation of isomorphism coincides with the orbit equivalence relation of
the canonical action of the group $S_{\infty }$ of permutations of $\mathbb{N%
}$. Conversely the orbit equivalence relation associated with a continuous
action of $S_{\infty }$ on a Polish space can be regarded as the relation of
isomorphism within some Borel class of $\mathcal{L}$-structures for a
suitable language $\mathcal{L}$ \cite[Theorem 3.6.1]{gao_invariant_2009}.

The natural next step beyond countable structures is to consider orbit
equivalence relations of arbitrary Polish group actions. An equivalence
relation is then \emph{classifiable by orbits }if it is Borel reducible to
the orbit equivalence relation associated with a Polish group action. In
view of the L\'{o}pez-Escobar theorem for continuous logic \cite%
{coskey_lopez-escobar_2014,ben_yaacov_lopez-escobar_2014}, any orbit
equivalence relation of a Polish group action can be identified with the
relation of isomorphism within some Borel class of separable $\mathcal{L}$%
-structures for a language $\mathcal{L}$ in the logic for metric structures 
\cite{ben_yaacov_model_2008, farah_model_?}. Conversely any such a relation
of isomorphism is Borel reducible to the action of a Polish group \cite%
{elliott_isomorphism_2013}. Beside classifiability by countable structures
and by orbits, there are many other natural benchmarks of complexity, such
as those corresponding to the classes of analytic equivalence relations that
are hyperfinite, hypersmooth, treeable, or countable.

The notion of Borel reducibility leads to a very rich and interesting
theory, an account of which can be found in \cite{gao_invariant_2009}.
Albeit highly mixed and interlaced, one can isolate two main directions of
research. One concerns the general study of the hierarchy of Borel (or
analytic) equivalence relations under Borel reducibility. Remarkable
examples of this approach are the many deep dichotomy results for Borel
equivalence relations, starting from the seminal results of Silver \cite%
{silver_counting_1980} and Harrington-Kechris-Louveau \cite%
{harrington_glimm-effros_1990}. A survey on this line of research can be
found in \cite{hjorth_new_1997,kechris_new_1999,miller_graph-theoretic_2012}.

The other direction aims at the study of important examples of equivalence
relations coming from the practice of mathematics in various fields. The
goal of this investigation is to determine as precisely as possible where
these concrete classification problems sit in the Borel complexity
hierarchy. Such an analysis is valuable in that it provides information on
which classification results are possible, and what is the right kind of
invariants to classify a given class of objects.

Beyond the smooth case, many explicit classification results in mathematics
involve countable structures up to isomorphism as invariants. It is
therefore of great importance to understand when an equivalence relation is
classifiable by countable structure. Suppose that the equivalence relation
under consideration is the orbit equivalence relation $E_{G}^{X}$ of a
Polish group action $G\curvearrowright X$.\ The above description of
isomorphism relations of countable structures as actions of $S_{\infty }$
suggest that one can use dynamical properties of the action $%
G\curvearrowright X$ to rule out the existence of a Borel reduction of $%
E_{G}^{X}$ to $E_{Y}^{S_{\infty }}$ for any $S_{\infty }$-space $Y$. For
example it is well known that if there is a dense orbit, and all orbits are
meager (a condition known as \emph{generic ergodicity}) then the orbit
equivalence relation $E_{G}^{X}$ is not smooth by \cite[Proposition 6.1.10]%
{gao_invariant_2009}. By locally strengthening this condition, Hjorth
introduced the notion of \emph{(generically) turbulent }action.

Suppose that $X$ is a Polish $G$-space, $x\in X$, $U$ is an open
neighborhood of $x$ in $X$, and $V$ is an open neighborhood of the identity
in $G$. The \emph{local orbit} $\mathcal{O}\left( x,U,V\right) $ is the set
of all points that can be reached from $x$ by applying elements of $V$
without ever leaving $U$. A point $x$ of $X$ is \emph{turbulent }if for
every nonempty open subset $U$ of $X$ and every open neighborhood $V$ of the
identity of $G$ the \emph{local orbit }$\mathcal{O}\left( x,U,V\right) $ is
somewhere dense. The $G$-orbit of a turbulent points is a \emph{turbulent
orbit}.

\begin{definition}[Hjorth, 2000]
The action $G\curvearrowright X$ is \emph{turbulent }if every orbit is
dense, meager, and turbulent. It is \emph{generically turbulent }if there is
a $G$-invariant $G_{\delta }$ subset $C$ of $X$ such that the restriction of
the action of $G$ to $C$ is turbulent.
\end{definition}

Since the set of turbulent points is a $G$-invariant $G_{\delta }$ set, a $G$%
-space is generically ergodic if and only if every orbit is meager, and
there is a dense turbulent orbit. The following is the main result of
Hjorth's turbulence theory.

\begin{theorem}[Hjorth, 2000]
\label{Theorem:Hjorth-turbulence}Suppose that $G\curvearrowright X$ is a
generically turbulent Polish group action. If $Y$ is an $S_{\infty }$-space,
and $f:X\rightarrow Y$ is a Baire measurable function that maps $G$-orbits
into $S_{\infty }$-orbits, then there is an $S_{\infty }$-orbit whose
preimage under $f$ is comeager. In particular $E_{G}^{X}$ is not
classifiable by countable structures.
\end{theorem}

Hjorth's theory of turbulence is currently the main tool to establish
results of nonclassifiability by countable structures. A more general
approach in terms of forcing has been recently developed by Zapletal \cite%
{zapletal_forcing_2014, zapletal_analytic_2013}.

Among the classifications problems that have been proved intractable by
means of turbulence, there are unitary equivalence of irreducible
representations of locally compact groups \cite{hjorth_non-smooth_1997} and
C*-algebras \cite{farah_dichotomy_2012, kerr_turbulence_2010}, conjugacy of
ergodic measure-preserving transformations of the standard probability space 
\cite{hjorth_invariants_2001,foreman_anti-classification_2004}, and unitary
conjugacy of self-adjoint and unitary operators on the Hilbert space \cite%
{kechris_strong_2001}.

All the equivalence relations mentioned above are classifiable by the orbits
of Polish group actions. Many naturally occurring classification problems
have been shown to transcend such a benchmark, and in fact being of maximum
complexity for all analytic equivalence relations. These relations that are 
\emph{complete }for analytic equivalence relations include isomorphism and
(complemented) biembeddability of Banach spaces, isomorphism of Polish
groups, uniform homeomorphism of complete separable metric spaces \cite%
{ferenczi_complexity_2009}, isomorphism of operator spaces \cite%
{argerami_classification_2014}, bi-embeddability of countable graphs \cite%
{louveau_complete_2005}.

Considering the importance that classification---in the form of the Elliott
classification program \cite%
{elliott_classification_1995,elliott_regularity_2008}---has in the modern
theory of C*-algebras and C*-dynamics, it is not surprising that Borel
complexity theory has provided interesting new information on the subject.
The descriptive analysis of the classification problem for C*-algebra has
been initiated in \cite%
{camerlo_completeness_2001,farah_turbulence_2014,farah_descriptive_2012},
and continued in \cite%
{elliott_isomorphism_2013,sabok_completeness_2013,zielinski_complexity_2014}%
. In this survey we will focus on some new results on the Borel complexity
of \emph{automorphisms }of C*-algebras, mostly from the papers \cite%
{lupini_unitary_2014,kerr_borel_2014,gardella_conjugacy_2014}.

\subsection*{Acknowledgments}

We would like to thank Ilijas Farah for his several comments on earlier
versions of this survey that---we believe---significantly contributed to
improve the presentation of the material. We are also grateful to Samuel Coskey for sharing his notes on the relation of isomorphism of countable torsion-free abelian groups.

\section{C*-algebras and noncommutative dynamics\label{Section:introduction}}

An (abstract) \emph{C*-algebra} $A$ is a complex algebra endowed with a norm
and an involution $x\mapsto x^{\ast }$ satisfying $\left\Vert xy\right\Vert
\leq \left\Vert x\right\Vert \left\Vert y\right\Vert $ and the C*-identity $%
\left\Vert x^{\ast }x\right\Vert =\left\Vert x\right\Vert ^{2}$ for $x,y\in
A $. Suppose that $H$ is a Hilbert space. Let $B\left( H\right) $ be the
algebra of bounded linear operators on $H$. Considering on $B\left( H\right) 
$ the operator norm and the involution that assigns to every operator its
adjoint\emph{\ }yields a C*-algebra structure on $B\left( H\right) $. A
(concrete) C*-algebra is then a closed subalgebra of $B\left( H\right) $
that is self-adjoint, i.e.\ closed under taking adjoints. The foundational
result of the theory of C*-algebra due to Gelfand and Naimark asserts that
any abstract C*-algebra is isomorphic to a concrete one \cite%
{gelfand_imbedding_1943}.

For simplicity \emph{we will assume all C*-algebras to be unital}, i.e.\ to
contain a multiplicative identity $1$, and \emph{separable}.\emph{\ }These
can be concretely identified with the separable closed self-adjoint
subalgebras of $B\left( H\right) $ that contain the identity operator. There
are many good references for the theory of C*-algebras \cite%
{blackadar_operator_2006,brown_c*-algebras_2008,davidson_c*-algebras_1996}.
In the following we will mostly refer to \cite{blackadar_operator_2006} for
convenience.

The original motivation in the development of the theory of C*-algebras was
to provide rigorous mathematical foundation to quantum mechanics \cite%
{segal_c*-algebras_1994}. In essence the \textquotedblleft quantization
process\textquotedblright\ consists in replacing functions with operators on
the Hilbert space. If $K$ is a compact metrizable space, then the space $%
C\left( K\right) $ of complex-valued continuous functions on $K$ is a
C*-algebra with respect to the pointwise operations and the uniform norm.
The C*-algebras of this form are precisely the \emph{abelian }ones \cite[%
II.2.2.4]{blackadar_operator_2006}. Arbitrary C*-algebras can therefore be
thought as quantized or noncommutative spaces. The theory of C*-algebras is
thus sometimes referred to as \emph{noncommutative topology}. The fact that
these quantized spaces have no actual points is consistent with the
principle that it is physically meaningless to speak of a \textquotedblleft
point\textquotedblright\ in the phase space of a quantum particle \cite%
{effros_some_1994}.

A similar quantization process can be applied to \emph{topological dynamical
systems}. Classically a topological dynamical system is given by a pair $%
\left( K,T\right) $ where $K$ is a compact metrizable space and $T$ is a
homeomorphism of $K$. Dynamical systems can be thought as symmetries of a
space, but also as specifying the (discrete) time evolution of a system. In
order to understand what is the right quantized analog of dynamical systems,
one can observe that a homeomorphism $T$ of $K$ canonically induces a map $%
\alpha _{T}:C\left( K\right) \rightarrow C\left( K\right) $ defined by $%
f\mapsto f\circ T$. Such a map is an \emph{automorphism }of $C\left(
K\right) $, in the sense that it preserves all the C*-algebra structure of $%
C\left( K\right) $. Conversely given any automorphism $\alpha $ of $C\left(
K\right) $ there is a unique homeomorphism $T$ of $K$ that corresponds to $%
\alpha $ via the above relation. Therefore the quantized version of a
topological dynamical system is a pair $\left( A,\alpha \right) $ where $A$
is a C*-algebra, and $\alpha $ is an automorphism of $A$. The study of
automorphisms of C*-algebras can also be described as\emph{\ noncommutative
topological dynamics}.

The automorphisms of a C*-algebra $A$ form a group $\mathrm{Aut}\left(
A\right) $ under composition. Such a group is endowed with a canonical
Polish group topology, which is the weakest topology making the
point-evaluations continuous. This observation opens up the possibility to
the use of descriptive set-theoretic methods in the study of automorphisms
of C*-algebras.

In the classical setting, topological dynamical systems are naturally
studied and classified up to the relation of \emph{conjugacy}. Two systems $%
T,T^{\prime }$ on the same space $K$ are conjugate if they are equal modulo
a symmetry of the space. This means that there is an invertible continuous
transformation $S$ of $K$ such that $S\circ T=T^{\prime }\circ S$. There is
a straightforward quantized version of such a notion: two automorphisms $%
\alpha ,\alpha ^{\prime }$ of a C*-algebra $A$ are conjugate if there is an
automorphism $\gamma $ of $A$ such that $\alpha \circ \gamma =\gamma \circ
\alpha ^{\prime }$. Equivalently, $\alpha $ and $\alpha ^{\prime }$ are
conjugate inside the automorphism group $\mathrm{Aut}\left( A\right) $ of $A$%
. However this is not the only possible analog of conjugacy in the
noncommutative setting. In order to define the other one, we need to
introduce the important notion of unitary element.

An element $u$ of $A$ is \emph{unitary }if it satisfies $uu^{\ast }=u^{\ast
}u=1$. When $A$ is concretely represented as an algebras of operators on $H$%
, this is equivalent to the assertion that $u$ is a unitary operator on $H$,
i.e.\ a surjective linear isometry. The unitary elements of $A$ form a
multiplicative group $U\left( A\right) $ which is Polish when endowed with
the norm topology. Any unitary element $u$ induces an automorphism $\mathrm{%
Ad}\left( u\right) $ of $A$ defined by $x\mapsto uxu^{\ast }$. Automorphisms
of this form are called \emph{inner }and form a normal subgroup $\mathrm{Inn}%
\left( A\right) $ of $\mathrm{Aut}\left( A\right) $. These are the
automorphisms of $A$ that come from \textquotedblleft
rotating\textquotedblright\ the Hilbert space $H$ one which $A$ is
represented, and can be regarded as \emph{trivial}. It is therefore natural,
in trying to classify the automorphisms of $A$, not to distinguish between
automorphisms that are equal up to an inner automorphism. This leads to the
notion of unitary equivalence of automorphisms of $A$.

Two automorphisms $\alpha ,\beta $ of $A$ are \emph{unitarily equivalent} if 
$\alpha =\mathrm{Ad}\left( u\right) \circ \beta $ for some unitary $u$ of $A$%
. Equivalently the images of $\alpha $ and $\beta $ in the quotient group $%
\mathrm{Aut}\left( A\right) /\,\mathrm{Inn}\left( A\right) $ are equal.
Combining unitary equivalence with conjugacy leads to the definition of
cocycle conjugacy. Two automorphisms are \emph{cocycle conjugate} if one is
unitarily equivalent to a conjugate of the other. Equivalently their images
in the quotient group $\mathrm{Aut}\left( A\right) /\,\mathrm{Inn}\left(
A\right) $ are conjugate. The name cocycle conjugacy comes from the fact
that it can be described in terms of perturbations by a cocycle.

In the commutative setting, conjugacy and cocycle conjugacy coincide since
there are no nontrivial inner automorphisms. In the noncommutative setting
they behave in general very differently. It could be matter of debate which
one is the \textquotedblleft right\textquotedblright\ quantized analog of
conjugacy of homeomorphisms. It is an empiric fact that most satisfactory
classification results concern the relation of cocycle conjugacy. This is
not a coincidence, as it has recently been confirmed in \cite%
{kerr_borel_2014} by means of Hjorth's theory of turbulence that the
relation of conjugacy is not well behaved for \textquotedblleft
most\textquotedblright\ C*-algebras.

Similar conclusions hold for the relation of unitary equivalence. In fact it
has recently been shown in \cite{lupini_unitary_2014} that such an
equivalence relation is not classifiable by countable structures whenever it
is not smooth. Such a dichotomy is reminiscent of the analogous phenomenon
for the relation of unitary equivalence of irreducible representations of
C*-algebras. A classical result of Glimm asserts that such a relation is
smooth precisely when $A$ is type I \cite{glimm_type_1961}. Recently Glimm's
result has been refined by showing that in the non type I case the
irreducible representations are not classifiable by countable structures 
\cite{hjorth_non-smooth_1997,kerr_turbulence_2010,farah_dichotomy_2012}.

The situation for the relation of cocycle conjugacy if far less clear. Many
positive classification results have been obtained in special cases
including, most notably, automorphisms of Kirchberg algebras satisfying
suitable freeness conditions \cite{nakamura_aperiodic_2000}. It has recently
been shown in \cite{gardella_conjugacy_2014} that the relation of cocycle
conjugacy of automorphisms of the Cuntz algebra $\mathcal{O}_{2}$ is a
complete analytic set. This is obtained by analyzing the range of the
invariant in a classification result of Izumi from \cite{izumi_finite_2004}
and proving the existence of a Borel inverse for the classifying function.
Such a result provides one of the few known lower bounds on the complexity
of the classification problem for automorphisms of C*-algebras up to cocycle
conjugacy.

\section{The Elliott classification program\label{Section:nuclear}}

\subsection{Tensor products\label{Subsection:tensor}}

Tensor products play a key role in the theory of C*-algebras. Suppose that $%
A $ and $B$ are C*-algebras. There is a natural way to define a norm on the
algebraic tensor product $A\odot B$ of $A$ and $B$ as complex algebras:
concretely represent $A$ and $B$ as closed self-adjoint subalgebras of $%
B\left( H\right) $. Consider the inclusion of $A\odot B$ into $B\left(
H\otimes H\right) $ given by the action $\left( a\otimes b\right) \left( \xi
\otimes \eta \right) =a\xi \otimes a\eta $ for $a\in A$, $b\in B$, and $\xi
,\eta \in H$. Here $H\overline{\otimes }H$ denotes the Hilbert space tensor
product of $H$ by itself, i.e.\ the completion of the algebraic tensor
product $H\odot H$ with respect to the pre-inner product%
\begin{equation*}
\left\langle \xi \otimes \eta ,\xi ^{\prime }\otimes \eta ^{\prime
}\right\rangle =\left\langle \xi ,\xi ^{\prime }\right\rangle \left\langle
\eta ,\eta ^{\prime }\right\rangle \text{.}
\end{equation*}%
The inclusion $A\odot B\subset B\left( H\right) $ yields a cross C*-norm on $%
A\odot B$, i.e.\ a norm such that $\left\Vert a\otimes b\right\Vert
=\left\Vert a\right\Vert \left\Vert b\right\Vert $ and the completion of $%
A\odot B$ with respect to such a norm is a C*-algebra. In general there may
many other cross C*-norms on $A\odot B$. A C*-algebra $A$ for which this
does \emph{not }happen is called nuclear. More precisely a C*-algebra $A$ is 
\emph{nuclear }if for any other C*-algebra $B$, there is a unique cross
C*-norm on $A\odot B$. The corresponding completion is then simply denoted
by $A\otimes B$.

All abelian C*-algebras are nuclear, and in fact the unique tensor product
is given by the formula%
\begin{equation*}
C\left( X\right) \otimes C\left( Y\right) \cong C\left( X\times Y\right)
\end{equation*}%
where $X\times Y$ denotes the product space. Similarly, full matrix algebras 
$M_{n}\left( \mathbb{C}\right) $ are nuclear, and the tensor product $%
M_{n}\left( \mathbb{C}\right) \otimes M_{k}\left( \mathbb{C}\right) $ can be
identified with $M_{nk}\left( \mathbb{C}\right) $.

The standard example of a C*-algebra that is not nuclear is the reduced
C*-algebra $C_{r}^{\ast }\left( \mathbb{F}_{2}\right) $ of the free group on 
$2$ generators. This is the C*-algebra of operators on the Hilbert space $%
\ell ^{2}\left( \mathbb{F}_{2}\right) $---with canonical basis $\left( \xi
_{h}\right) _{h\in \mathbb{F}_{2}}$---generated by the unitary operators $%
u_{g}\left( \xi _{h}\right) =\xi _{gh}$ for $g\in \mathbb{F}_{2}$. In fact
the tensor product $C_{r}^{\ast }\left( \mathbb{F}_{2}\right) \odot
C_{r}^{\ast }\left( \mathbb{F}_{2}\right) $ has the largest possible number
of cross C*-norms: continuum many \cite{wiersma_c*-norms_2014}. Another
important example of C*-algebra that is not nuclear is $B\left( H\right) $
when $H$ is the separable infinite-dimensional Hilbert space \cite%
{wassermann_tensor_1976}. Again the amount of cross C*-norms on $B\left(
H\right) \odot B\left( H\right) $ as large as possible: power of the
continuum \cite{ozawa_continuum_2014}.

Nuclearity is a key notion in the theory of C*-algebras. Many equivalent
reformulations of nuclearity have been given including Banach-algebraic
amenability, the completely positive approximation property, and injectivity
or semidiscreteness of the second dual \cite[IV.3.1]{blackadar_operator_2006}%
. In the following we will only consider tensor products of nuclear
C*-algebras. A nuclear C*-algebra $A$ is $B$-\emph{absorbing} if $A\otimes
B\cong A$ and \emph{self-absorbing} if $A\otimes A\cong A$.

The Elliott classification program is an ambitious project aiming at a
complete classification of (algebraically) simple nuclear C*-algebras by
K-theoretic data \cite{elliott_classification_1995, elliott_regularity_2008}.

\subsection{The Elliott invariant\label{Subsection:Elliott}}

K-theory was first developed in algebraic geometry by Atiyah and Hirzebruch 
\cite{atiyah_riemann-roch_1959} to study vector bundles by algebraic means
(topological\emph{\ }K-theory).\ These ideas were then translated into an
algebraic language, leading to algebraic\emph{\ }K-theory of Banach algebras
and general rings \cite{rosenberg_algebraic_1994}. The K-theoretic machinery
was then incorporated in the theory of C*-algebras regarded as noncommutative%
\emph{\ }topology \cite{schochet_algebraic_1994}.

The K-theory of a C*-algebra $A$ can be described in terms of the $K_{0}$
and $K_{1}$ groups.\ The ordered $K_{0}$-group $K_{0}\left( A\right) $ can
be obtained from the \emph{Murray--von Neumann semigroup} $V\left( A\right) $
as the group of formal differences, regarding $V\left( A\right) $ as the
cone of positive elements inside of $K_{0}\left( A\right) $. Such a
construction is entirely analogous to the construction of the ordered group $%
\mathbb{Z}$ from the additive semigroup $\mathbb{N}$. The semigroup $V\left(
A\right) $ can be described---algebraically---as the set of equivalence
classes of finitely generated projective modules over $A$ up to stable
equivalence, where the operation corresponds to the direct sum of modules.
Equivalently one can describe $V\left( A\right) $ intrinsically in terms of
the relation of Murray-von Neumann equivalence of projections. (All
projections are assumed to be orthogonal projections onto their range. They
can be characterized algebraically as those operators satisfying $p=p^{\ast
}=p^{2}$.) Two projections $p,q\in A$ are equivalent in $A$ if $A$ contains
a \emph{partial isometry }$v$ with support projection $p$ and range
projection $q$. This means that $v$ is an isometry from the range of $p$
onto the range of $q$, and identically zero on the orthogonal complement.
Equivalently, $v^{\ast }v=p$ and $vv^{\ast }=q$. Such a notion captures the
fact that $p$ and $q$ have the same dimension \emph{relatively to }$A$. It
was first considered by Murray and von\ Neumann in their classification
theory for von Neumann factors \cite{murray_rings_1943}. In order to define $%
V\left( A\right) $ one needs to consider all the matrix amplifications $%
M_{n}\left( A\right) $ and their projections up to equivalence. The
semigroup operation correspond to \textquotedblleft putting projections on
the diagonal\textquotedblright , in formulas%
\begin{equation*}
p\oplus q=%
\begin{bmatrix}
p & 0 \\ 
0 & q%
\end{bmatrix}%
\text{.}
\end{equation*}%
When $A\cong C\left( X\right) $ is abelian, projections on matrix
amplifications of $A$ correspond to \emph{vector bundles }over $X$. The
notion of Murray-von Neumann equivalence of projections in this case
translates into the natural notion of stable equivalence of vector bundles.
The $K_{0}$-group thus codes the higher rank topological information on $A$.
The $K_{1}$-group $K_{1}\left( A\right) $ is defined similarly by
considering, instead of projections, unitary elements up to the relation of
being connected by a continuous path. More information on K-theory of
C*-algebras can be found in \cite{blackadar_K-theory_1986}.

The \emph{Elliott invariant }of a C*-algebra $A$ comprises the $K_{0}$ and $%
K_{1}$ groups and the \emph{trace simplex }$T\left( A\right) $ of $A$. A 
\emph{trace }on $A$ is a linear functional $\phi $ on $A$ of norm $1$
satisfying $\phi \left( 1\right) =1$ and $\phi \left( ab\right) =\phi \left(
ba\right) $ for $a,b\in A$. This can be seen as a noncommutative analog of a
Radon probability measure. In fact when $A\cong C\left( X\right) $ is
abelian, any Radon probability measure $\mu $ on $X$ yields a trace on $%
C\left( X\right) $ given by the integral functional $\int \left( \cdot
\right) d\mu $, and any trace is of this form. The trace simplex $T\left(
A\right) $ is a Choquet simplex with respect to the topological and affine
structure obtained by regarding of $T\left( A\right) $ as a convex subset of
the dual $A^{\ast }$ of $A$.

The last bit of information included in the Elliott invariant is the
canonical pairing between $K_{0}$ and traces, obtained by evaluating a trace
at (any representative of) an equivalence class of projections. Such a
pairing records how the measure-theoretic and the topological structures of
the algebra interact. The Elliott invariant is furthermore \emph{functorial}%
: a morphism between algebras induces a morphism at the level of invariants.

\subsection{AF and UHF algebras\label{Subsection:AF}}

One of the first instances of classification of a large important class of
C*-algebras by K-theoretic invariants is the Elliott-Bratteli classification
of approximately finite-dimensional (AF) algebras \cite%
{bratteli_inductive_1972,elliott_classification_1976}. For algebras in this
class it is shown that any morphisms between the invariants can be \emph{%
lifted }to a morphism between the algebras. The main proof technique, which
plays an important role in many other further classification results, is an 
\emph{approximate intertwining argument }that can be regarded as the
C*-algebraic version of Cantor's back-and-forth method.

The uniformly hyperfinite (UHF) algebras are a particularly important class
of AF algebras. These are infinite tensor products of algebras of the form $%
M_{n}\left( \mathbb{C}\right) $. This class of algebras had been previously
classified by Glimm \label{glimm_certain_1960} in terms of the corresponding 
\emph{generalized natural number}. Such a number has \textquotedblleft prime
factors\textquotedblright\ $p^{n}$ where $n\in \mathbb{N}\cup \left\{ \infty
\right\} $ is the supremum of those $k$ such that $M_{p^{k}}(\mathbb{C})$
unitally embeds into the algebra. A UHF algebra is self-absorbing if and
only if it is of \emph{infinite type}, i.e.\ every prime appears in the
associated generalized natural number with infinite exponent.

\subsection{Kirchberg algebras\label{Subsection:Kirchberg}}

The Elliott-Bratteli classification result for AF algebras was later
extended to many other classes of algebras. Most notably a complete
classification of purely infinite\emph{\ }simple nuclear C*-algebras (\emph{%
Kirchberg algebras}) was obtained by Kirchberg and Phillips in \cite%
{kirchberg_exact_1995, kirchberg_embedding_2000} modulo the technical
assumption that they satisfy the Universal Coefficient Theorem (UCT). It
should be mentioned that the UCT holds in all the currently known examples
of simple nuclear purely infinite C*-algebras. In fact it is currently a
major open problem whether there is any nuclear C*-algebra for which the UCT
fails.

The notion of purely infinite C*-algebra is defined in terms of the notion
of \emph{infinite projection}. A projection $p$ in a C*-algebra $A$ is \emph{%
infinite }if it is Murray-von Neumann equivalent to a \emph{proper
subprojection} of $p$. This is a projection $q$ distinct from $p$ and
dominated by $p$ in the sense that the range of $q$ is contained in the
range of $p$ or, equivalently, $pq=qp=q$. Such a definition should be
compared with the fact that a set is infinite if and only if it is in
bijection with a proper subset (Dedekind-infiniteness). A simple C*-algebra
is then \emph{purely infinite} if it has an abundance of infinite
projections: for every nonzero element $x$ of $A$ the subalgebra $x^{\ast }Ax
$ contains an infinite projection. A \emph{Kirchberg algebra }is a simple
nuclear purely infinite C*-algebra. Kirchberg algebras do not have any trace %
\label{rordam_stable_2004}. Moreover the order structure on the $K_{0}$%
-group trivializes, and every element is positive. Therefore the Elliott
invariant of a Kirchberg algebra $A$ reduces to the pair of (countable,
discrete) abelian groups $K_{0}\left( A\right) $ and $K_{1}\left( A\right) $.

Fundamental examples of purely infinite simple nuclear C*-algebras are the 
\emph{Cuntz algebras }$\mathcal{O}_{2}$ and $\mathcal{O}_{\infty }$. The
Cuntz algebra $\mathcal{O}_{2}$ is the C*-algebra generated by two
isometries $v_{1},v_{2}$ with complementary ranges. This means that $%
v_{1},v_{2}\in B\left( H\right) $ are linear isometries with images two
complementary subspaces of $H$ or, equivalently, $v_{1}v_{1}^{\ast
}+v_{2}v_{2}^{\ast }=v_{1}^{\ast }v_{1}=v_{2}^{\ast }v_{2}=1$. The Cuntz
algebra $\mathcal{O}_{\infty }$ is generated by infinitely many isometries $%
\left( v_{n}\right) _{n\in \mathbb{N}}$ with mutually orthogonal ranges,
i.e.\ satisfying $v_{i}^{\ast }v_{i}=1$ and $v_{i}^{\ast }v_{j}=0$ for $%
i\neq j$. One can similarly define the Cuntz algebra $\mathcal{O}_{n}$ for
any $n\geq 2$. Such algebras were first defined and studied by Cuntz in \cite%
{cuntz_simple_1977}. Cuntz showed that the isomorphism class of these
algebras does not depend on the choice of the generating isometries.
Moreover they are simple and purely infinite. Their K-theory was computed in 
\cite{cuntz_K-theory_1981}: they all have trivial $K_{1}$-group, while $%
K_{0}\left( \mathcal{O}_{\infty }\right) \cong \mathbb{Z}$ and $K_{0}\left( 
\mathcal{O}_{2}\right) $ is trivial. More generally $K_{0}\left( \mathcal{O}%
_{n}\right) \cong \mathbb{Z}/\!\left( n-1\right) \mathbb{Z}$. The first step
in the Kirchberg-Phillips classification consists in showing that, if $B$ is
a simple nuclear C*-algebra, then $\mathcal{O}_{2}$ is $B$-absorbing, and $B$
is purely infinite if and only if it is $\mathcal{O}_{\infty }$-absorbing.

\subsection{The Jiang-Su algebra\label{Subsection:Jiang-Su}}

The possible extent of the Elliott classification program has been recently
limited by breakthroughs due R\o rdam and Toms. In \cite%
{rordam_simple_2003,toms_independence_2005, toms_classification_2008} they
showed that in general the Elliott invariant is \emph{not }a complete
invariant for simple nuclear C*-algebras. In fact continuum many pairwise
nonisomorphic simple nuclear C*-algebras with the same Elliott invariant are
constructed in \cite{toms_comparison_2009}. These algebras are distinguished
by another invariant, the \emph{radius of comparison}, which is a
noncommutative analog of the \emph{mean dimension }of dynamical systems \cite%
{lindenstrauss_mean_1999,lindenstrauss_mean_2000}. The radius of comparison
is detected by the first order theory of the C*-algebra seen as a structure
in continuous logic \cite{ben_yaacov_model_2008, farah_model_?}. It
therefore remains open the---unlikely---possibility that the Elliott
invariant \emph{together with the first order theory }provides a complete
invariant for simple nuclear C*-algebras; see \cite{farah_model_?-1}.

Many efforts have been recently dedicated to isolate the class of
\textquotedblleft well behaved\textquotedblright\ C*-algebras for which the
Elliott classification program can be successfully recasted. In these
efforts the \emph{Jiang-Su algebra }$\mathcal{Z}$ has a role of paramount
importance. This is an infinite-dimensional simple nuclear C*-algebra that
is a tensorial \textquotedblleft neutral element\textquotedblright\ at the
level of invariants. It was first constructed by Jiang and Su in \cite%
{jiang_simple_1999} as limit of \emph{dimension drop algebras}. These are
algebras of the form%
\begin{equation*}
\left\{ f\in C\left( \left[ 0,1\right] ,M_{p}\left( \mathbb{C}\right)
\otimes M_{q}\left( \mathbb{C}\right) \right) :f\left( 0\right) \in
M_{p}\left( \mathbb{C}\right) \otimes 1,f\left( 1\right) \in 1\otimes
M_{q}\left( \mathbb{C}\right) \right\}
\end{equation*}%
endowed with the canonical trace coming from the Lebesgue measure on $\left[
0,1\right] $. The Jiang-Su construction was shown to be of Fra\"{\i}ss\'{e}%
-theoretical nature in \cite{eagle_fraisse_2014}.

If $A$ is any simple nuclear C*-algebra, then $A\otimes \mathcal{Z}$ and $A$
have the same Elliott invariant. It follows that it is a \emph{necessary
condition} for $A$ to be classifiable that $A$ be $\mathcal{Z}$-\emph{%
absorbing}. The notion of $\mathcal{Z}$-absorption is the subject of an
important conjecture due to Toms and Winter. Such a conjecture asserts that,
for simple nuclear C*-algebras, $\mathcal{Z}$-absorption is equivalent to
other two regularity properties of topological and cohomological nature,
respectively: finite nuclear dimension and strict comparison. Such a
conjecture has been by now verified in a vast number of cases including all
C*-algebras whose trace simplex is \textquotedblleft small\textquotedblright
; see \cite{bosa_covering_2014} and references therein.

\subsection{Strongly self-absorbing C*-algebras\label{Subsection:ssa}}

The Jiang-Su algebra $\mathcal{Z}$, the Cuntz algebras $\mathcal{O}_{2}$ and 
$\mathcal{O}_{\infty }$, and the infinite-type UHF\ algebras all play a key
role in the classification program of C*-algebras. The fundamental
properties that make these algebras stand out have been isolated in \cite%
{toms_strongly_2007}. A (necessarily nuclear) C*-algebra is \emph{strongly
self-absorbing} if the embedding $a\mapsto a\otimes 1$ of $A$ into $A\otimes
A$ is approximately unitarily equivalent to an isomorphism. This means that
there exist an isomorphism $\gamma :A\rightarrow A\otimes A$ and a sequence $%
\left( u_{n}\right) $ of unitaries of $A$ such that $\left\Vert \gamma
\left( a\right) -u_{n}\left( a\otimes 1\right) u_{n}^{\ast }\right\Vert
\rightarrow 0$. A strongly self-absorbing C*-algebra $A$ is in particular
self-absorbing, and in fact isomorphic to the infinite tensor product $%
A^{\otimes \mathbb{N}}$ of copies of $A$. The algebras $\mathcal{Z}$, $%
\mathcal{O}_{2}$, $\mathcal{O}_{\infty }$, the infinite type UHF algebras
and their tensor products are the only currently known examples of strongly
self-absorbing C*-algebras.

\subsection{Borel complexity of C*-algebras\label{Subsection:Borel}}

The analysis of the classification problem of C*-algebras from the
perspective of Borel complexity theory has been initiated in \cite%
{farah_turbulence_2014}. The Elliott invariant has been shown to be \emph{%
Borel-computable} with respect to any natural parametrization of C*-algebras 
\cite{farah_descriptive_2012}. As observed in \cite[\S 3]{farah_logic_2014},
by combining results from \cite{elliott_isomorphism_2013,
melleray_computing_2007, gao_classification_2003, sabok_completeness_2013,
zielinski_complexity_2014} one can conclude that the following relations all
have maximum complexity among equivalence relations that are classifiable by
the orbits of a Polish group action:

\begin{itemize}
\item Isomorphism of (Elliott-classifiable) C*-algebras;

\item Homeomorphism of compact metrizable spaces;

\item Affine homeomorphism of Choquet simplices;

\item Isometry of Banach spaces;

\item Isometry of complete separable metric spaces.
\end{itemize}

While some of the classification problems above are perceived to be less
intractable than others, they have in fact the same complexity from the
perspective of invariant descriptive set theory. However the functorial
nature of the classification problem for C*-algebras is not taken into
account by this analysis. \emph{Functorial Borel complexity} is a refinement
introduced in \cite{lupini_polish_2014} of the usual notion of Borel
complexity, aiming at capturing the complexity of classifying the objects of
a category in an explicit and functorial way. It is conceivable that some of
the classes above might have distinct functorial complexity.

\section{Conjugacy\label{Section:conjugacy}}

In this section we will consider the relation of conjugacy for automorphisms
of a C*-algebra $A$. When $A$ is abelian, $A\cong C\left( K\right) $ for
some compact metrizable space, and $\mathrm{Aut}\left( A\right) \cong 
\mathrm{Homeo}\left( K\right) $. Therefore the relation of conjugacy for
automorphisms of $A$ coincides with the relation of conjugacy of
homeomorphisms of $K$. The latter equivalence relation has been studied for
particular instances of $K$. For example it has been shown in \cite[\S 4.2]%
{hjorth_classification_2000} that the homeomorphisms of the unit interval $%
\left[ 0,1\right] $ are classifiable by countable structures. However an
increase of dimension entails an increase of complexity: the homeomorphisms
of the unit square $\left[ 0,1\right] ^{2}$ are not classifiable by
countable structures \cite[4.3]{hjorth_classification_2000}. This result was
one of the first applications of Hjorth's turbulence theory, developed
therein.

When $K$ is \emph{zero-dimensional}, $K$ is the Stone space of a countable
Boolean algebra $\mathcal{B}$, and $\mathrm{Homeo}\left( K\right) \cong 
\mathrm{Aut}\left( \mathcal{B}\right) $. Therefore conjugacy of
homeomorphisms of $K$ coincides with conjugacy of automorphisms of $\mathcal{%
B}$. Since $\mathrm{Aut}\left( \mathcal{B}\right) $ is a closed subgroup of $%
S_{\infty }$, it follows that the latter equivalence relation is
classifiable by countable structures. When $\mathcal{B}$ is the countable
atomless Boolean algebra, its Stone space $K$ is the Cantor set. The
relation of conjugacy of homeomorphisms of $K$ has been shown to have
maximum complexity among relations that are classifiable by countable
structures in \cite{camerlo_completeness_2001}.

Some positive classification results up to conjugacy for automorphisms of
finite order by K-theoretic invariants have been obtained by Izumi \cite%
{izumi_finite_2004,izumi_finite_2004-1}. These results assume the
automorphisms (or their \textquotedblleft dual actions\textquotedblright )
to have suitable freeness conditions such as the Rokhlin property \cite%
{phillips_tracial_2012}. Moreover they apply to classes of algebras that are
classifiable\ in the sense of the Elliott program, such as the class of
Kirchberg algebras; see \S \ref{Subsection:Kirchberg}. The starting point
of\ Izumi analysis is the Kirchberg-Phillips classification result, and the
study of automorphisms at the level of the invariants.

In the rest of this section we will present the following nonclassifiability
result for automorphisms of $\mathcal{Z}$-absorbing C*-algebras.

\begin{theorem}[Kerr-Lupini-Phillips, 2014]
\label{Theorem:conjugacy}Suppose that $A$ is a C*-algebra such that $%
A\otimes \mathcal{Z}\cong A$. Then the automorphisms of $A$ are not
classifiable up to conjugacy by countable structures.
\end{theorem}

Theorem \ref{Theorem:conjugacy} in particular rules out classification by
K-theoretic data in view of the Borel-computability of K-theory \cite%
{farah_descriptive_2012}. As observed in \S \ref{Subsection:Jiang-Su}, the
class of $\mathcal{Z}$-absorbing C*-algebras contains all C*-algebras that
fall under the scope of the classification program. The proof strategy for
Theorem \ref{Theorem:conjugacy} is to establish generic turbulence of the
conjugation action $\mathrm{Aut}\left( \mathcal{Z}\right) \curvearrowright 
\mathrm{Aut}\left( \mathcal{Z}\right) $, and then use Hjorth's turbulence
theorem together with the canonical embedding $\mathrm{Aut}\left( \mathcal{Z}%
\right) \hookrightarrow \mathrm{Aut}\left( A\otimes \mathcal{Z}\right) \cong 
\mathrm{Aut}\left( A\right) $ to deduce nonclassifiability in $\mathrm{Aut}%
\left( A\right) $.

Recall from \S \ref{Subsection:ssa} that the Jiang-Su algebra is a
particular instance of an important class of algebras known as strongly
self-absorbing \cite[Definition 1.1]{toms_strongly_2007}. Such a property in
particular implies that $\mathcal{Z}$ is isomorphic to the infinite tensor
product $\mathcal{Z}^{\otimes \mathbb{Z}}$ of copies of $\mathcal{Z}$. As a
consequence after replacing $\mathcal{Z}$ with $\mathcal{Z}^{\otimes \mathbb{%
Z}}$ one can consider the \emph{shift} automorphism $\sigma $ of $\mathcal{Z}%
^{\otimes \mathbb{N}}$ associated with the shift $n\mapsto n-1$ on $\mathbb{Z%
}$. A key property of the shift $\sigma $ is being \emph{malleable }\cite[%
Lemma 2.6]{kerr_borel_2014}. This means that there is a continuous path $%
\left( \rho _{t}\right) _{t\in \left[ 0,1\right] }$ in $\mathrm{Aut}\left( 
\mathcal{Z}^{\otimes \mathbb{Z}}\otimes \mathcal{Z}^{\otimes \mathbb{Z}%
}\right) $ such that $\rho _{0}$ is the identity, $\rho _{1}$ is the flip $%
x\otimes y\mapsto y\otimes x$, and $\rho _{t}\circ \left( \sigma \otimes
\sigma \right) =\left( \sigma \otimes \sigma \right) \circ \rho _{t}$ for
every $t\in \left[ 0,1\right] $ \cite[Definition 2.5]{kerr_borel_2014}. The
notion of malleability plays a key role in Popa's deformation-rigidity
theory for von Neumann algebras and II$_{1}$ factors \cite%
{popa_deformation_2007}. Such a notion turns out to be of fundamental
importance also in this C*-algebraic setting. In particular it is the key
property that allows one to establish the following result.

\begin{proposition}
\label{Proposition:shift}The the shift automorphism $\sigma $ of $\mathcal{Z}%
^{\otimes \mathbb{Z}}\cong \mathcal{Z}$ is a turbulent point with dense
orbit for the conjugation action $\mathrm{Aut}\left( \mathcal{Z}\right)
\curvearrowright \mathrm{Aut}\left( \mathcal{Z}\right) $.
\end{proposition}

Proposition \ref{Proposition:shift} is obtained via a factor exchange
argument applied to the tensor product $\mathcal{Z}^{\otimes \mathbb{Z}%
}\otimes \mathcal{Z}^{\otimes \mathbb{Z}}$. Malleability of the shift is
used to carry out the exchange via a continuous path of unitaries in a way
that commutes with the shift itself. Density of the shift is a consequence
of work of Sato on \emph{stability}\ properties of the shift and, more
generally, automorphisms of $\mathcal{Z}$ satisfying a suitable freeness
condition \cite{sato_rohlin_2010}.

In order to conclude the proof that the conjugation action $\mathrm{Aut}%
\left( \mathcal{Z}\right) \curvearrowright \mathrm{Aut}\left( \mathcal{Z}%
\right) $ is generically turbulent, one then needs show that every orbit is
meager. This can be established as an application of a result of Rosendal
that infers meagerness of conjugacy classes from a density condition for
periodic elements.

\begin{lemma}[Rosendal]
\label{Lemma:Rosendal}Suppose that $G$ is a Polish group. If for every
infinite subset $I$ of $\mathbb{N}$ the set of $g\in G$ such that $g^{n}=1$
for some $n\in I$ is dense, then $G$ has meager conjugacy classes.
\end{lemma}

The hypothesis of Lemma \ref{Lemma:Rosendal} imply the existence of a
collection of conjugation-invariant comeager subsets of $G$ with trivial
intersection. This directly imply meagerness of conjugacy classes; see \cite[%
Proposition 18]{rosendal_generic_2009} and \cite[page 9]{kechris_global_2010}%
.

In order to conclude that $\mathrm{\mathrm{Aut}}\left( \mathcal{Z}\right) $
satisfies the hypothesis of Lemma \ref{Lemma:Rosendal}, one can observe that
the shift $\sigma $ can be approximated by \textquotedblleft periodic
shifts\textquotedblright\ of sufficiently large period. Formally these are
the automorphisms of $\mathcal{Z}^{\otimes \mathbb{N}}\cong \mathcal{Z}$
associated with cyclic permutations of $\mathbb{Z}$ 
\begin{equation*}
\bigl(%
\begin{smallmatrix}
-n & -n+1 & -n+2 & \cdots  & m-1 & m%
\end{smallmatrix}%
\bigr)
\end{equation*}%
for $n,m\in \mathbb{N}$. This suffices in view of the fact that the shift
has dense conjugacy class.

We now explain how to deduce Theorem \ref{Theorem:conjugacy} from
Proposition \ref{Proposition:shift}. Suppose that $A$ is a $\mathcal{Z}$%
-absorbing C*-algebra. Identifying $A$ with $A\otimes \mathcal{Z}$ we can
regard $\mathrm{Aut}\left( \mathcal{Z}\right) $ as a closed subgroup of $%
\mathrm{Aut}\left( A\otimes \mathcal{Z}\right) $ via the inclusion $\gamma
\mapsto id_{A}\otimes \gamma $. One can then argue by contradiction that, if
the automorphisms of $A\otimes \mathcal{Z}$ are classifiable by countable
structures, then the relation for elements of $\mathrm{Aut}\left( \mathcal{Z}%
\right) $ of being conjugate inside $\mathrm{Aut}\left( A\otimes \mathcal{Z}%
\right) $ must have a comeager class. However the latter fact can be
excluded essentially via the same argument that shows that $\mathrm{Aut}%
\left( \mathcal{Z}\right) $ has meager conjugacy classes.

\section{Unitary equivalence\label{Section:unitary-equivalence}}

In this section we consider the relation of unitary equivalence of
automorphisms of C*-algebras. Suppose that $A$ is a C*-algebra. Recall that
any unitary element $u$ of $A$ induces an automorphism $\mathrm{Ad}\left(
u\right) $ of $A$ defined by $x\mapsto uxu^{\ast }$. Automorphisms of this
form are called inner and form a normal subgroup $\mathrm{Inn}\left(
A\right) $ of $\mathrm{Aut}\left( A\right) $. Two automorphisms of $A$ are 
\emph{unitarily equivalent }if they belong to the same coset of $\mathrm{Inn}%
\left( A\right) $. It should be noted that $\mathrm{Inn}\left( A\right) $ is
the continuous homomorphic image of the unitary group $U\left( A\right) $ of 
$A$ under the map $u\mapsto \mathrm{Ad}\left( u\right) $. Therefore, being
the continuous homomorphic image of a Polish group, it is a Borel subset of $%
\mathrm{Aut}\left( A\right) $ \cite[Exercise 15.15]{kechris_classical_1995}.
Moreover a standard argument in descriptive set theory allows one to
conclude that the relation of unitary equivalence of automorphisms of $A$ is 
\emph{smooth }if and only if $\mathrm{Inn}\left( A\right) $ is a \emph{closed%
} subgroup of $\mathrm{Aut}\left( A\right) $. In fact if $\mathrm{Inn}\left(
A\right) $ is closed then the quotient group $\mathrm{Aut}\left( A\right) /\,%
\mathrm{Inn}\left( A\right) $ is Polish, and its elements parametrize the
cosets of $\mathrm{Inn}\left( A\right) $. Conversely if unitary equivalence
is smooth then by turning Borel sets into clopen sets \cite[\S 13.A]%
{kechris_classical_1995} one can find a finer topology on $\mathrm{Aut}%
\left( A\right) $ that makes such $\mathrm{Inn}\left( A\right) $ closed. By
the uniqueness of the Polish topology in a Polish group, such a topology
coincides with the original one, and hence $\mathrm{Inn}\left( A\right) $
was already closed.

An explicit classification of automorphisms of $A$ up to unitary equivalence
has been obtained in \cite{phillips_automorphisms_1980} when $A$ has
continuous trace. In particular this classification shows that in the
continuous trace case such a relation is smooth. Conversely it has been
shown in \cite{phillips_outer_1987} that if $A$ does not have continuous
trace, then the relation of unitary equivalence of automorphisms is not
smooth. The last result has been strengthened in \cite{lupini_unitary_2014}
by showing that when $A$ does not have continuous trace, then in fact the
automorphisms of $A$ are not classifiable up to unitary equivalence by
countable structures. In this section we will present an overview of these
results, after introducing the notion of continuous trace C*-algebra.

A \emph{central sequence }in a C*-algebra $A$ is a bounded sequence $\left(
a_{n}\right) $ in $A$ such that the commutators $a_{n}b-ba_{n}$ converge to $%
0$ for any $b\in A$. A central sequence is \emph{trivial }if there is a
sequence $\left( z_{n}\right) $ in the center of $A$ such that $\left\Vert
a_{n}-z_{n}\right\Vert \rightarrow 0$ when $n\rightarrow +\infty $. (The 
\emph{center} $Z\left( A\right) $ of $A$ is the subalgebra of those elements
of $A$ that commute with any other element of $A$.) The C*-algebra $A$ \emph{%
has continuous trace} if every central sequence in $A$ is trivial. While
this is not the original definition of continuous trace C*-algebras, it has
been proved to be equivalent by Akemann and Pedersen \cite%
{akemann_central_1979} building on previous work of Elliott \cite%
{elliott_some_1977}.

The fact that $\mathrm{Inn}\left( A\right) $ is closed in $\mathrm{Aut}%
\left( A\right) $ if and only if $A$ has continuous trace---observed in \cite%
{elliott_some_1977, phillips_outer_1987}---follows from a standard fact in
the theory of C*-algebras: a C*-algebra is generated by its unitary elements 
\cite[II.3.2.16]{blackadar_operator_2006}. In particular the center of its
unitary group coincides with the intersection of $U\left( A\right) $ with
the center of $A$. Recall also the classical principles from descriptive set
theory that a Borel bijection between standard Borel spaces is a Borel
isomorphism \cite[Corollary 15.2]{kechris_classical_1995}, and a Borel
homomorphism between Polish groups is continuous \cite[Exercise 9.16]%
{kechris_classical_1995}.

If $\mathrm{Inn}\left( A\right) $ is closed in $\mathrm{Aut}\left( A\right) $%
, then the map $w\mapsto \mathrm{\mathrm{Ad}}\left( w\right) $ from the
unitary group $U\left( A\right) $ of $A$ to $\mathrm{Inn}\left( A\right) $
induces, passing to the quotient, a Polish group isomorphism from the
quotient of $U\left( A\right) $ by its center to $\mathrm{Inn}\left(
A\right) $. If $\left( w_{n}\right) $ is a central sequence in $A$, then $%
\mathrm{Ad}\left( w_{n}\right) $ converges to the identity automorphism of $A
$. Henceforth $\left( w_{n}\right) $ converges to the identity in the
quotient of $U\left( A\right) $ by its center. This means that there is a
sequence $\left( z_{n}\right) $ in $Z\left( A\right) \cap U\left( A\right) $
such that $\left\Vert u_{n}-z_{n}\right\Vert \rightarrow 0$. In other words $%
\left( w_{n}\right) $ is a trivial central sequence. The converse
implication can be proved just reversing the argument above.

Such a result admits the following strengthening, which can be seen as an
instance of the principle that \textquotedblleft when it is bad, it is
worse\textquotedblright ; see \cite[Theorem 1.1]{lupini_unitary_2014}.

\begin{theorem}
\label{Theorem:nonclassification-ue}If $A$ is a C*-algebra that does not
have continuous trace, then the automorphisms of $A$ are not classifiable up
to unitary equivalence by countable structures.
\end{theorem}

Theorem \ref{Theorem:nonclassification-ue} in particular implies a dichotomy
for the complexity of the relation of unitary equivalence of automorphisms
of C*-algebras. Such a relation is either smooth, or not classifiable by
countable structures. The same phenomenon has been observed for the relation
of unitary equivalence of irreducible representations of locally compact
groups and C*-algebras \cite{hjorth_non-smooth_1997, farah_dichotomy_2012,
kerr_turbulence_2010}.

In the rest of this section we will present the main ideas in the proof of
Theorem \ref{Theorem:nonclassification-ue}. Again Hjorth's theory of
turbulence will play a central role. One of the first examples of turbulent
action provided by Hjorth is the action of $\ell ^{1}$ on $\mathbb{R}^{%
\mathbb{N}}$ by translation \cite[\S 3.3]{hjorth_classification_2000}. Here $%
\ell ^{1}$ is regarded as an (additive) Polish group endowed with the
topology induced by the $\ell ^{1}$-norm, while $\mathbb{R}^{\mathbb{N}}$ is
endowed with the product topology.

Hjorth's turbulence theorem together with turbulence of the translation
action of $\ell ^{1}$ on $\mathbb{R}^{\mathbb{N}}$ yields the following
nonclassifiability criterion; see \cite[Criterion 3.3]{lupini_unitary_2014}.

\begin{criterion}
\label{Criterion:nonclassifiability}Suppose that $E$ is an equivalence
relation on a standard Borel space $X$.\ If there is a Borel function $%
f:\left( 0,1\right) ^{\mathbb{N}}\rightarrow X$ such that

\begin{enumerate}
\item $f\left( \mathbf{x}\right) E{}f\left( \mathbf{x}^{\prime }\right) $
whenever $\mathbf{x},\mathbf{x}^{\prime }\in \left( 0,1\right) ^{\mathbb{N}}$
are such that $\mathbf{x}-\mathbf{x}^{\prime }\in \ell ^{1}$, and

\item the preimage under $f$ of any $E$-class is meager,
\end{enumerate}

then $E$ is not classifiable by countable structures.
\end{criterion}

Criterion \ref{Criterion:nonclassifiability} will be the main tool to
establish Theorem \ref{Theorem:nonclassification-ue}. Let $A$ be a
C*-algebra that does not have continuous trace. Thus $A$ contains a
nontrivial central sequence $\left( x_{n}\right) $. We will assume for
simplicity that $\left( x_{n}\right) $ is nontrivial in the following
stronger sense: there is a central sequence $\left( y_{n}\right) $ in $A$
such that the commutators $x_{n}y_{n}-y_{n}x_{n}$ do not converge to $0$.
This covers most interesting cases, including when the C*-algebra $A$ is
simple.

After replacing $x_{n}$ with $x_{n}^{\ast }x_{n}$ we can assume that the $%
x_{n}$'s are \emph{positive}. This means that they are positive linear
operators in any concrete representation of $A$ as a subalgebra of $B\left(
H\right) $. The idea is now to use the sequence $\left( x_{n}\right) $ to
produce, with an eye at Criterion \ref{Criterion:nonclassifiability},
\textquotedblleft many\textquotedblright\ automorphisms of $A$ indexed by
sequences $\left( t_{n}\right) $ in $\left( 0,1\right) ^{\mathbb{N}}$. In
order to define these automorphisms we need to turn the $x_{n}$'s into
unitary elements. This can be achieved by taking exponentials. Recall that
if $x$ is a positive real number then $\exp \left( ix\right) =\sum_{k\in 
\mathbb{N}}\frac{1}{k!}\left( ix\right) ^{k}$ is a complex number of modulus 
$1$, i.e.\ a unitary element of $\mathbb{C}$.\ The same fact---with similar
proof---applies when $\mathbb{C}$ is replaced by an arbitrary C*-algebra. If 
$x$ is a positive element of $A$, then $\exp \left( ix\right) =\sum_{k\in 
\mathbb{N}}\frac{1}{k!}\left( ix^{k}\right) $ is a well defined unitary
element of $A$.

We can then define for $\left( t_{n}\right) \in \left( 0,1\right) ^{\mathbb{N%
}}$ an automorphism $\alpha _{\left( t_{n}\right) }$ of $A$ as the limit of
the sequence of inner automorphisms $\mathrm{Ad}\left( \exp \left(
it_{n}x_{n}\right) \right) $ associated with the unitary elements $\exp
\left( it_{n}x_{n}\right) $. The fact that $\left( x_{n}\right) $ is a
central sequence guarantees, perhaps after passing to a subsequence, that
such automorphisms are well defined. The map $f:\mathbf{t}\mapsto \alpha _{%
\mathbf{t}}$ is easily seen to be Borel (in fact, continuous). Moreover if $%
\mathbf{t},\mathbf{s}\in \left( 0,1\right) ^{\mathbb{N}}$ are such that $%
\mathbf{t}-\mathbf{s}\in \ell ^{1}$ then $\alpha _{\mathbf{t}}=\mathrm{Ad}%
\left( u\right) \circ \alpha _{\mathbf{s}}$ where $u$ is the limit inside
the unitary group of $A$ of the sequence%
\begin{equation*}
\left( \exp \left( it_{1}x_{1}\right) \cdots \exp \left( it_{n}x_{n}\right)
\exp \left( -is_{n}x_{n}\right) \cdots \exp \left( -is_{1}x_{1}\right)
\right) _{n\in \mathbb{N}}\text{.}
\end{equation*}%
Thus $f$ satisfies Condition (1) of Criterion \ref%
{Criterion:nonclassifiability}. In order to establish Condition (2), one has
to show that the preimage under $f$ of any unitary equivalence class is
meager. If $X\subset \left( 0,1\right) ^{\mathbb{N}}$ is nonmeager, then it
is easy to see that there are $\mathbf{t},\mathbf{s}\in X$ such that $%
\left\vert t_{n}-s_{n}\right\vert >\frac{1}{2}$ for infinitely many $n\in 
\mathbb{N}$. One can then use the sequence $\left( y_{n}\right) $ to argue
that, for such a choice of $\mathbf{s}$ and $\mathbf{t}$ in $X$, $\alpha _{%
\mathbf{s}}\circ \alpha _{\mathbf{t}}^{-1}$ is not inner. (More details can
be found in \cite[Proposition 5.6]{lupini_unitary_2014}.) Therefore $\alpha
_{\mathbf{s}}$ and $\alpha _{\mathbf{t}}$ are not unitarily equivalent. This
concludes the proof that the preimage under $f$ of any unitary equivalence
class is meager, and $f$ satisfies the assumptions of Criterion \ref%
{Criterion:nonclassifiability}.

\section{Cocycle conjugacy\label{Section:cocycle-conjugacy}}

Recall that two automorphisms $\alpha $ and $\beta $ are cocycle conjugate
if a conjugate of $\alpha $ is unitarily equivalent to $\beta $. In formulas 
$\alpha =\mathrm{Ad}\left( u\right) \circ \gamma \circ \beta \circ \gamma
^{-1}$ for some unitary element $u$ of $A$ and some automorphism $\gamma $
of $A$. The majority of positive classification results for automorphisms of
C*-algebras are up to cocycle conjugacy. Yet, not much is known in general,
and most of the results hold under the additional assumption that the
automorphisms have the (weak) Rokhlin property. Such a property prescribes
the existence of a noncommutative analog of the Rokhlin tower from ergodic
theory. It turns out that under very general assumptions automorphisms with
the Rokhlin property exist, and in fact form a dense $G_{\delta }$ set \cite%
{phillips_tracial_2012}. Little is known about the complexity of the
relation of cocycle conjugacy outside the (comeager) set of automorphism
with the Rokhlin property.

Kishimoto has shown in \cite{kishimoto_rohlin_1995} that for UHF algebras
any two automorphisms with the Rokhlin property are cocycle conjugate. A
similar result has been obtained by Sato for automorphisms of the Jiang-Su
algebra \cite{sato_rohlin_2010}. In both cases the proof relies on the
characterization of automorphisms with the Rokhlin property as automorphisms
all of whose powers are outer in a strong sense. This indicates that the
Rokhlin property can be seen as a quantized analog of freeness for
topological dynamical systems. The other fundamental feature of automorphism
with the Rokhlin property that is used in Kishimoto and Sato's result is 
\emph{stability}. Stability guarantees that if $\alpha $ is an automorphism
with the Rokhlin property and $u$ is a unitary, then $u$ can be approximated
with a unitary of the form $v\alpha \left( v^{\ast }\right) $ for some
unitary $v$. This allows one two perturb any two automorphisms with the
Rokhlin property to be cocycle conjugate.

In the case of Kirchberg algebras satisfying the Universal Coefficient
Theorem, the Kirchberg-Phillips classification mentioned in \S \ref%
{Subsection:Kirchberg} also provides information about automorphisms. Recall
that such a classification result asserts that UCT Kirchberg algebras are
classified by their Elliott invariants. Since Kirchberg algebras have no
traces, the Elliott invariant reduces to the $K_{0}$ and $K_{1}$ groups. A
key ingredient of the Kirchberg-Phillips classification result consists in
showing that two morphisms between Kirchberg algebras induces the same
morphism at the level of invariants if and only if they are \emph{%
asymptotically }unitarily equivalent. Such a notion is the analog of unitary
equivalence where one replaces a single unitary with a continuous path. It
was later shown by Nakamura that, under the additional assumption of the
Rokhlin property, asymptotic unitary equivalence implies cocycle conjugacy 
\cite[Theorem 5]{nakamura_aperiodic_2000}. Again the stability of Rokhlin
property automorphisms plays a major role. Combining Nakamura's result and
the Kirchberg-Phillips classification, one obtains that Rokhlin property
automorphisms of UCT Kirchberg algebras are cocycle conjugate if and only if
the corresponding morphisms at the level of the Elliott invariant are
conjugate \cite[Theorem 9]{nakamura_aperiodic_2000}.

In the particular case of the Cuntz algebra $\mathcal{O}_{2}$, both the $%
K_{0}$ and the $K_{1}$ groups are trivial. Therefore it follows from
Nakamura's classification that any two automorphisms of $\mathcal{O}_{2}$
with the Rokhlin property are cocycle conjugate. Beyond the class of
automorphisms with the Rokhlin property, a classification result has been
obtained by Izumi for strongly approximately inner symmetries of $\mathcal{O}%
_{2}$ \cite{izumi_finite_2004, izumi_finite_2004-1}. Again the
Kirchberg-Phillips classification of UCT Kirchberg algebras plays a key
role, together with Kirchberg's absorption theorem asserting that $\mathcal{O%
}_{2}\otimes B\cong \mathcal{O}_{2}$ for any simple nuclear C*-algebra $B$.

A \emph{symmetry} is an automorphism of order $2$. A symmetry is \emph{%
strongly approximately inner }if it is approximated by inner automorphisms
induced by unitaries that are approximately of order $2$ and are
approximately fixed by the given automorphism. It should be mentioned that
no symmetry of $\mathcal{O}_{2}$ that is not strongly approximately inner is
currently known. In fact very recent results suggest that the problem of the
existence of such symmetries is tightly connected with the problem of the
existence of nuclear C*-algebras for which the UCT fails.

The invariant involved in Izumi's classification is the\emph{\ crossed
product}, which is a C*-algebra encoding the action. The crossed product is
one of the most important constructions in C*-dynamics \cite[\S II.10]%
{blackadar_operator_2006}. Succinctly, if $A$ is a C*-algebra and $\alpha $
is an automorphism, then the corresponding crossed product $A\rtimes
_{\alpha }\mathbb{Z}$ is the \textquotedblleft universal\textquotedblright\
C*-algebra generated by $A$ and a unitary $u$ that \textquotedblleft
implements $\alpha $\textquotedblright\ in the sense that $uxu^{\ast
}=\alpha \left( x\right) $ for every $x\in A$.

It follows from the Izumi and Kirchberg-Phillips classification results
that, under the additional assumption that the crossed product satisfies the
UCT, the K-theory of the crossed product provides a complete invariant up to
cocycle conjugacy for strongly approximately inner symmetries of $\mathcal{O}%
_{2}$. The range of the invariants has also been computed by Izumi. This is
the collection of all pairs of \emph{countable uniquely }$2$\emph{-divisible
abelian groups}. (An abelian group is $2$-divisible if every element can be
written as double of some element. When such a representation is unique, the
group is uniquely $2$-divisible.) The fact that such a classifying map is
Borel can be seen using the methods from \cite{farah_descriptive_2012}. It
has been shown in \cite{gardella_conjugacy_2014} that such a map admits a
Borel inverse. For simplicity the attention is restricted there to algebras
with trivial $K_{1}$-group.

\begin{proposition}[Gardella-Lupini, 2014]
There is a Borel map $G\mapsto \alpha _{G}$ from the space of uniquely $2$%
-divisible countable abelian group to the space of strongly approximately
inner symmetries of $\mathcal{O}_{2}$ such that $G$ is the $K_{0}$-group of
the crossed product of $\mathcal{O}_{2}$ by $\alpha _{G}$, and the $K_{1}$%
-group of such a crossed product is trivial.
\end{proposition}

Such a result together with Izumi's classification shows the following:

\begin{theorem}[Gardella-Lupini, 2014]
\label{Theorem:reduction-symmetry}The relation of isomorphism of uniquely $2$%
-divisible countable abelian groups is Borel reducible to the relation of
cocycle conjugacy of (strongly approximately inner) symmetries of $\mathcal{O%
}_{2}$.
\end{theorem}

A minor modification of the argument Downey-Montalb\'{a}n from \cite%
{downey_isomorphism_2008}---see also \cite{hjorth_isomorphism_2002}---shows
that the relation of isomorphism of torsion free $2$-divisible countable
abelian group is a complete analytic set.

Recall that if $A$, $B$ are analytic subsets of standard Borel spaces $X$, $%
Y $, then $A$ is \emph{Wadge-reducible} to $B$, in formulas $A\leq _{W}B$,
if there is a Borel function $f:X\rightarrow Y$ such that $A$ is the
preimage of $B$ under $f$. An analytic set is \emph{complete }if it is a $%
\leq _{W}$-maximum element of the class of analytic sets. The well known
fact that there are analytic sets that are not Borel shows in particular
that a complete analytic set is not Borel.

The canonical example of a complete analytic set is the set of \emph{%
ill-founded trees }on $\mathbb{N}$ inside the space of all trees. (A tree is
ill-founded if it has an infinite branch, and well-founded otherwise.) In
order to show that the relation of isomorphism of $2$-divisible torsion-free
abelian groups is a complete analytic set, it is enough to find a Wadge
reduction from the set of ill-founded trees to such a relation. To
facilitate the task one can use the following theorem: there are Borel maps $%
R,S$ from trees to trees such that, for any tree $T$

\begin{itemize}
\item if $T$ is ill-founded, then $R\left( T\right) \cong S\left( T\right) $
is ill-founded, and

\item if $T$ is well-founded, then $R\left( T\right) $ is well-founded and $%
S\left( T\right) $ is ill-founded.
\end{itemize}

In fact it suffices to define maps $R,S$ as above from trees to \emph{%
countable linear orders}. One can then just take the tree of finite
descending chains of the linear order. The existence of such maps has been
observed in \cite[Lemma 3.1]{downey_isomorphism_2008}. The key idea is to
use Harrison's theorem from \cite{harrison_recursive_1968}, asserting that a 
\emph{pseudo-well order} $L$, i.e.\ a countable linear with no infinite
descending chains that are hyperarithmetic in $L$, is ill-founded if and
only if it is of the form $\omega _{1}^{L}\left( 1+\mathbb{Q}\right) +\alpha 
$ for some $\alpha <\omega _{1}^{L}$. To define $R\left( T\right) $, first
pair up $T$ with the complete binary tree. This gives a tree $T_{1}$ with
either none or uncountably many infinite branches. Then express the set of
branches of $T_{1}$ that are not hyperarithmetic in $T$ as $\left\{ x\in
\omega ^{<\omega }:\widehat{T}\left( x\right) \text{ is ill-founded}\right\} 
$ for some tree $\widehat{T}$ on $\omega \times \omega $ computable in $T$.
(Here $\widehat{T}\left( x\right) $ denotes the section tree of $\widehat{T}$
corresponding to $x$.) Let now $L$ be the Brower-Kleene order of $\widehat{T}
$, and $R\left( T\right) =L+L+L+\cdots =L\omega $. Observe that if $T$ is
well-founded, then $R\left( T\right) $ is well-founded, otherwise $R\left(
T\right) \cong \omega _{1}^{T}\left( 1+\mathbb{Q}\right) $ by Harrison's
theorem. Similarly one can define a Borel-computable map such that $S\left(
T\right) \cong \omega _{1}^{T}\left( 1+\mathbb{Q}\right) $ for every $T$.

In view of such a result, to conclude the proof that the isomorphism
relation of countable $2$-divisible torsion-free abelian groups is a
complete analytic set, it is enough to assign in a Borel fashion to any tree 
$T$ a countable $2$-divisible torsion free abelian group $T_{G}$ in such a
way that

\begin{itemize}
\item if $T\cong T^{\prime }$ then $G_{T}\cong G_{T^{\prime }}$, and

\item if $G_{T}\cong G_{T^{\prime }}$ then exactly one between $T$ and $%
T^{\prime }$ is well-founded.
\end{itemize}

To perform such a construction, fix enumerations $\left( p_{n}\right) $ and $%
\left( q_{n}\right) $ of pairwise disjoint sets of odd primes. Suppose that $%
T$ is a tree, with vertex set $V$. Let $\mathbb{Q}^{\left( V\right) }$ be
the group of formal sums $\sum_{v\in V}q_{v}v$ where $q_{v}\in \mathbb{Q}$
and all but finitely many $q_{v}$'s are zero. Define $G_{T}$ to be the
subgroup of $\mathbb{Q}^{\left( V\right) }$ spanned by elements of the form $%
\frac{1}{2p_{n}^{k}}v$ and $\frac{1}{2q_{n}^{k}}\left( v+w\right) $ where $%
k\in \mathbb{N}$, $v$ is a vertex of $T$ of height $n$, and $w$ is an
immediate successor of $v$ in $T$.

It is clear from the construction that isomorphic trees yield isomorphic
groups. In order to verify the second condition above, it is enough to argue
that one can detect the existence of an infinite branch of $T$ by looking at
group-theoretic properties of\ $G_{T}$. In fact one can verify that infinite
branches in $T$ correspond to sequences $\left( g_{n}\right) $ in $G_{T}$
such that $g_{n}$ is infinitely $p_{n}$-divisible, and $g_{n}+g_{n+1}$ is
infinitely $q_{n}$-divisible. This concludes the proof that the relation of
isomorphism of $2$-divisible torsion free abelian group is a complete
analytic set. From this and Theorem \ref{Theorem:reduction-symmetry} one can
deduce the following corollary.

\begin{corollary}
The relation of cocycle conjugacy of automorphisms of $\mathcal{O}_{2}$ is a
complete analytic set.
\end{corollary}

In fact the same holds for the relation of conjugacy of automorphisms of $%
\mathcal{O}_{2}$, essentially by the same argument. Moreover one can replace
symmetries with automorphisms of order $p$, and $2$-divisible groups with $p$%
-divisible groups for any prime $p$ using a construction of Barlak and Szab%
\'{o} \cite{barlak_rokhlin_2014}. Not much else is known in general about cocycle conjugacy of automorphisms
of $\mathcal{O}_{2}$.

\begin{problem}
\label{Problem:O2}Is the relation of cocycle conjugacy of automorphisms of $%
\mathcal{O}_{2}$ classifiable by countable structures?
\end{problem}

In fact in Problem \ref{Problem:O2} one can replace $\mathcal{O}_{2}$ with
any other simple nuclear C*-algebra; see \cite{gardella_conjugacy_2014}.

\providecommand{\bysame}{\leavevmode\hbox to3em{\hrulefill}\thinspace}
\providecommand{\MR}{\relax\ifhmode\unskip\space\fi MR }
\providecommand{\MRhref}[2]{%
  \href{http://www.ams.org/mathscinet-getitem?mr=#1}{#2}
}
\providecommand{\href}[2]{#2}


\begin{thebibliography}{10}

\bibitem{akemann_central_1979}
Charles~A. Akemann and Gert~K. Pedersen, \emph{Central sequences and inner
  derivations of separable {C}*-algebras}, American Journal of Mathematics
  \textbf{101} (1979), no.~5, 1047--1061.

\bibitem{argerami_classification_2014}
Martín Argerami, Samuel Coskey, Mehrdad Kalantar, Matthew Kennedy, Martino
  Lupini, and Marcin Sabok, \emph{The classification problem for finitely
  generated operator systems and spaces}, {arXiv}:1411.0512 (2014).

\bibitem{atiyah_riemann-roch_1959}
Michael~F. Atiyah and Friedrich Hirzebruch, \emph{Riemann-{R}och theorems for
  differentiable manifolds}, Bulletin of the American Mathematical Society
  \textbf{65} (1959), no.~4, 276--281.

\bibitem{barlak_rokhlin_2014}
Sel{\c{c}}uk Barlak and Gab{\'{o}}r Szab{\'{o}}, \emph{Rokhlin actions of
  finite groups on {UHF}-absorbing {C}*-algebras}, {arXiv}:1403.7312 (2014).

\bibitem{ben_yaacov_model_2008}
Ita{\"i} Ben~Yaacov, Alexander Berenstein, C.~Ward Henson, and Alexander
  Usvyatsov, \emph{Model theory for metric structures}, Model theory with
  applications to algebra and analysis. Vol. 2, London Mathematical Society
  Lecture Note Series, vol. 350, Cambridge University Press, 2008,
  p.~315{\textendash}427.

\bibitem{ben_yaacov_lopez-escobar_2014}
Ita{\"{i}} Ben~Yaacov, Andre Nies, and Todor Tsankov, \emph{A {L}opez-{E}scobar
  theorem for continuous logic}, {arXiv}:1407.7102 (2014).

\bibitem{blackadar_K-theory_1986}
Bruce Blackadar, \emph{{$K$}-theory for operator algebras}, Mathematical
  Sciences Research Institute Publications, vol.~5, Springer-Verlag, New York,
  1986.

\bibitem{blackadar_operator_2006}
\bysame, \emph{Operator algebras}, Encyclopaedia of Mathematical Sciences, vol.
  122, Springer-Verlag, Berlin, 2006.

\bibitem{bosa_covering_2014}
Joan Bosa, Nathanial~P. Brown, Yasuhiko Sako, Aaron Tikuisis, Stuart White, and
  Wilhelm Winter, \emph{Covering dimension of {C}*-algebras and the
  classification of maps by traces}, in preparation.

\bibitem{bratteli_inductive_1972}
Ola Bratteli, \emph{Inductive limits of finite dimensional {C}*-algebras},
  Transactions of the American Mathematical Society \textbf{171} (1972),
  195--234.

\bibitem{brown_c*-algebras_2008}
Nathanial~P. Brown and Narutaka Ozawa, \emph{C*-algebras and finite-dimensional
  approximations}, Graduate Studies in Mathematics, vol.~88, American
  Mathematical Society, 2008.

\bibitem{camerlo_completeness_2001}
Riccardo Camerlo and Su~Gao, \emph{The completeness of the isomorphism relation
  for countable {B}oolean algebras}, Transactions of the American Mathematical
  Society \textbf{353} (2001), no.~2, 491--518.

\bibitem{coskey_lopez-escobar_2014}
Samuel Coskey and Martino Lupini, \emph{A {L}{\'{o}}pez-{E}scobar theorem for
  metric structures, and the topological {V}aught conjecture},
  {arXiv}:1405.2859 (2014).

\bibitem{cuntz_simple_1977}
Joachim Cuntz, \emph{Simple {C}*-algebras generated by isometries},
  Communications in Mathematical Physics \textbf{57} (1977), no.~2, 173---185.

\bibitem{cuntz_K-theory_1981}
\bysame, \emph{K-theory for certain {C}*-algebras. {II}}, Journal of Operator
  Theory \textbf{5} (1981), no.~1, 101{\textendash}108.

\bibitem{davidson_c*-algebras_1996}
Kenneth~R. Davidson, \emph{C*-algebras by example}, Fields Institute
  Monographs, vol.~6, American Mathematical Society, 1996.

\bibitem{downey_isomorphism_2008}
Rod Downey and Antonio Montalb{\'a}n, \emph{The isomorphism problem for
  torsion-free abelian groups is analytic complete}, Journal of Algebra
  \textbf{320} (2008), no.~6, 2291--2300.

\bibitem{eagle_fraisse_2014}
Christopher~J. Eagle, Ilijas Farah, Bradd Hart, Boris Kadets, Vladyslav
  Kalashnyk, and Martino Lupini, \emph{Fra\"{i}ss\'{e} limits of
  {C}*-algebras}, {arXiv}:1411.4066 (2014).

\bibitem{effros_some_1994}
Edward~G. Effros, \emph{Some quantizations and reflections inspired by the
  {G}elfand-{N}a{\u{i}}mark theorem}, {C}*-algebras: 1943-1993 (San Antonio,
  {TX}, 1993), Contemporary Mathematics, vol. 167, Amer. Math. Soc.,
  Providence, {RI}, 1994, pp.~98--113.

\bibitem{elliott_classification_1976}
George~A Elliott, \emph{On the classification of inductive limits of sequences
  of semisimple finite-dimensional algebras}, Journal of Algebra \textbf{38}
  (1976), no.~1, 29--44.

\bibitem{elliott_some_1977}
George~A. Elliott, \emph{Some {C}*-algebras with outer derivations, {III}},
  Annals of Mathematics \textbf{106} (1977), no.~1, 121--143.

\bibitem{elliott_classification_1995}
\bysame, \emph{The classification problem for amenable {C}*-algebras},
  Proceedings of the International Congress of Mathematicians, Vol.\ 1, 2
  (Z{\"u}rich, 1994), Birkh{\"a}user, 1995, p.~922{\textendash}932.

\bibitem{elliott_isomorphism_2013}
George~A. Elliott, Ilijas Farah, Vern Paulsen, Christian Rosendal, Andrew~S.
  Toms, and Asger T{\"o}rnquist, \emph{The isomorphism relation for separable
  {C}*-algebras}, Mathematical Research Letters \textbf{20} (2013), no.~6,
  1071--1080.

\bibitem{elliott_regularity_2008}
George~A. Elliott and Andrew~S. Toms, \emph{Regularity properties in the
  classification program for separable amenable {C}*-algebras}, Bulletin of the
  American Mathematical Society \textbf{45} (2008), no.~2, 229--245.

\bibitem{farah_dichotomy_2012}
Ilijas Farah, \emph{A dichotomy for the {M}ackey {B}orel structure},
  Proceedings of the 11th Asian Logic Conference, World Sci. Publ., Hackensack,
  {NJ}, 2012, p.~86{\textendash}93.

\bibitem{farah_logic_2014}
\bysame, \emph{Logic and operator algebras}, Proceedings of the International
  Congress of Mathematicians (Seoul, South Corea), 2014.

\bibitem{farah_model_?-1}
Ilijas Farah, Bradd Hart, Martino Lupini, Leonel Robert, Aaron~P. Tikuisis,
  Alessandro Vignati, and Wilhelm Winter, \emph{Model theory of nuclear
  {C}*-algebras}, in preparation.

\bibitem{farah_model_?}
Ilijas Farah, Bradd Hart, and David Sherman, \emph{Model theory of operator
  algebras {II}: Model theory}, Israel Journal of Mathematics, to appear.

\bibitem{farah_descriptive_2012}
Ilijas Farah, Andrew~S. Toms, and Asger T{\"o}rnquist, \emph{The descriptive
  set theory of {C}*-algebra invariants}, International Mathematics Research
  Notices (2012), 5196--5226.

\bibitem{farah_turbulence_2014}
\bysame, \emph{Turbulence, orbit equivalence, and the classification of nuclear
  {C}*-algebras}, Journal f{\"u}r die reine und angewandte Mathematik
  \textbf{688} (2014), 101--146.

\bibitem{ferenczi_complexity_2009}
Valentin Ferenczi, Alain Louveau, and Christian Rosendal, \emph{The complexity
  of classifying separable banach spaces up to isomorphism}, Journal of the
  London Mathematical Society \textbf{79} (2009), no.~2, 323--345.

\bibitem{foreman_anti-classification_2004}
Matthew Foreman and Benjamin Weiss, \emph{An anti-classification theorem for
  ergodic measure preserving transformations}, Journal of the European
  Mathematical Society \textbf{6} (2004), no.~3, 277{\textendash}292.

\bibitem{friedman_borel_1989}
Harvey Friedman and Lee Stanley, \emph{A borel reductibility theory for classes
  of countable structures}, Journal of Symbolic Logic \textbf{54} (1989),
  no.~03, 894--914.

\bibitem{gao_invariant_2009}
Su~Gao, \emph{Invariant descriptive set theory}, Pure and Applied Mathematics
  (Boca Raton), vol. 293, {CRC} Press, Boca Raton, {FL}, 2009.

\bibitem{gao_classification_2003}
Su~Gao and Alexander~S. Kechris, \emph{On the classification of polish metric
  spaces up to isometry}, Memoirs of the American Mathematical Society
  \textbf{161} (2003), no.~766, 0--0.

\bibitem{gardella_conjugacy_2014}
Eusebio Gardella and Martino Lupini, \emph{Conjugacy and cocycle conjugacy of
  automorphisms of {$\mathcal{O}_2$} are not {B}orel}, {arXiv}:1404.3617
  (2014).

\bibitem{gelfand_imbedding_1943}
Izrail~M. Gelfand and Mark~A. Neumark, \emph{On the imbedding of normed rings
  into the ring of operators in {H}ilbert space}, Matematicheskii Sbornik
  Novaya Seriya \textbf{12(54)} (1943), 197{\textendash}213.

\bibitem{glimm_type_1961}
James Glimm, \emph{Type {I} {C}*-algebras}, Annals of Mathematics. Second
  Series \textbf{73} (1961), 572--612.

\bibitem{harrington_glimm-effros_1990}
Leo~A. Harrington, Alexander~S. Kechris, and Alain Louveau, \emph{A
  {G}limm-{E}ffros dichotomy for {B}orel equivalence relations}, Journal of the
  American Mathematical Society \textbf{3} (1990), no.~4, 903--928.

\bibitem{harrison_recursive_1968}
Joseph Harrison, \emph{Recursive pseudo-well-orderings}, Transactions of the
  American Mathematical Society \textbf{131} (1968), 526--543.

\bibitem{hjorth_non-smooth_1997}
Greg Hjorth, \emph{Non-smooth infinite-dimensional group representations},
  preprint, 1997.

\bibitem{hjorth_classification_2000}
\bysame, \emph{Classification and orbit equivalence relations}, Mathematical
  Surveys and Monographs, vol.~75, American Mathematical Society, Providence,
  {RI}, 2000.

\bibitem{hjorth_invariants_2001}
\bysame, \emph{On invariants for measure preserving transformations},
  Fundamenta Mathematicae \textbf{169} (2001), no.~1, 51{\textendash}84.

\bibitem{hjorth_isomorphism_2002}
\bysame, \emph{The isomorphism relation on countable torsion free abelian
  groups}, Fundamenta Mathematicae \textbf{175} (2002), no.~3,
  241{\textendash}257.

\bibitem{hjorth_new_1997}
Greg Hjorth and Alexander~S. Kechris, \emph{New dichotomies for {B}orel
  equivalence relations}, The Bulletin of Symbolic Logic \textbf{3} (1997),
  no.~3, 329--346.

\bibitem{izumi_finite_2004}
Masaki Izumi, \emph{Finite group actions on {C}*-algebras with the {R}ohlin
  property, {I}}, Duke Mathematical Journal \textbf{122} (2004), no.~2,
  233--280.

\bibitem{izumi_finite_2004-1}
\bysame, \emph{Finite group actions on {C}*-algebras with the {R}ohlin
  property{\textemdash}{II}}, Advances in Mathematics \textbf{184} (2004),
  no.~1, 119--160.

\bibitem{jiang_simple_1999}
Xinhui Jiang and Hongbing Su, \emph{On a simple unital projectionless
  {C}*-algebra}, American Journal of Mathematics \textbf{121} (1999), no.~2,
  359--413.

\bibitem{kechris_classical_1995}
Alexander~S. Kechris, \emph{Classical descriptive set theory}, Graduate Texts
  in Mathematics, vol. 156, Springer-Verlag, New York, 1995.

\bibitem{kechris_new_1999}
\bysame, \emph{New directions in descriptive set theory}, The Bulletin of
  Symbolic Logic \textbf{5} (1999), no.~2, 161--174.

\bibitem{kechris_global_2010}
\bysame, \emph{Global aspects of ergodic group actions}, Mathematical Surveys
  and Monographs, vol. 160, American Mathematical Society, Providence, {RI},
  2010.

\bibitem{kechris_strong_2001}
Alexander~S. Kechris and Nikolaos~E. Sofronidis, \emph{A strong generic
  ergodicity property of unitary and self-adjoint operators}, Ergodic Theory
  and Dynamical Systems \textbf{21} (2001), no.~5, 1459--1479.

\bibitem{kerr_turbulence_2010}
David Kerr, Hanfeng Li, and Mika{\"e}l Pichot, \emph{Turbulence,
  representations, and trace-preserving actions}, Proceedings of the London
  Mathematical Society \textbf{100} (2010), no.~2, 459--484.

\bibitem{kerr_borel_2014}
David Kerr, Martino Lupini, and N.~Christopher Phillips, \emph{Borel complexity
  and automorphisms of {C}*-algebras}, accepted for publication by the Journal
  of Functional Analysis.

\bibitem{kirchberg_exact_1995}
Eberhard Kirchberg, \emph{Exact {C}*-algebras, tensor products, and the
  classification of purely infinite algebras}, Proceedings of the International
  Congress of Mathematicians, Vol. 1, 2 (Z{\"u}rich, 1994), Birkh{\"a}user,
  Basel, 1995, p.~943{\textendash}954.

\bibitem{kirchberg_embedding_2000}
Eberhard Kirchberg and N.~Christopher Phillips, \emph{Embedding of exact
  {C}*-algebras in the {C}untz algebra {$\mathcal{O}_2$}}, Journal f{\"u}r die
  reine und angewandte Mathematik \textbf{525} (2000), 17--53.

\bibitem{kishimoto_rohlin_1995}
Akitaka Kishimoto, \emph{The {R}ohlin property for automorphisms of {UHF}
  algebras}, Journal f{\"{u}}r die reine und angewandte Mathematik
  \textbf{1995} (1995), no.~465, 183--196.

\bibitem{lindenstrauss_mean_1999}
Elon Lindenstrauss, \emph{Mean dimension, small entropy factors and an
  embedding theorem}, Institut des Hautes Études Scientifiques. Publications
  Mathématiques (1999), no.~89, 227--262 (2000).

\bibitem{lindenstrauss_mean_2000}
Elon Lindenstrauss and Benjamin Weiss, \emph{Mean topological dimension},
  Israel Journal of Mathematics \textbf{115} (2000), 1--24.

\bibitem{louveau_complete_2005}
Alain Louveau and Christian Rosendal, \emph{Complete analytic equivalence
  relations}, Transactions of the American Mathematical Society \textbf{357}
  (2005), no.~12, 4839--4866.

\bibitem{lupini_polish_2014}
Martino Lupini, \emph{Polish groupoids and functorial complexity},
  {arXiv}:1407.6671 (2014).

\bibitem{lupini_unitary_2014}
\bysame, \emph{Unitary equivalence of automorphisms of separable
  {C}*-algebras}, Advances in Mathematics \textbf{262} (2014), 1002--1034.

\bibitem{mackey_groups_1957}
George~W. Mackey, \emph{Borel structure in groups and their duals},
  Transactions of the American Mathematical Society \textbf{85} (1957), no.~1,
  134--165.

\bibitem{melleray_computing_2007}
Julien Melleray, \emph{Computing the complexity of the relation of isometry
  between separable banach spaces}, Mathematical Logic Quarterly \textbf{53}
  (2007), no.~2, 128--131.

\bibitem{miller_graph-theoretic_2012}
Benjamin~D. Miller, \emph{The graph-theoretic approach to descriptive set
  theory}, Bulletin of Symbolic Logic \textbf{18} (2012), no.~4, 554--575.

\bibitem{murray_rings_1943}
Francis~J. Murray and John von Neumann, \emph{On rings of operators. {IV}},
  Annals of Mathematics \textbf{44} (1943), no.~4, 716--808.

\bibitem{nakamura_aperiodic_2000}
Hideki Nakamura, \emph{Aperiodic automorphisms of nuclear purely infinite
  simple {C}*-algebras}, Ergodic Theory and Dynamical Systems \textbf{20}
  (2000), no.~6, 1749--1765.

\bibitem{ozawa_continuum_2014}
Narutaka Ozawa and Gilles Pisier, \emph{A continuum of c*-norms on
  {$\mathbb{B}(H)\otimes \mathbb{B}(H)$} and related tensor products}, to
  appear.

\bibitem{phillips_outer_1987}
John Phillips, \emph{Outer automorphisms of separable {C}*-algebras}, Journal
  of Functional Analysis \textbf{70} (1987), no.~1, 111--116.

\bibitem{phillips_automorphisms_1980}
John Phillips and Iain Raeburn, \emph{{Automorphisms of {C}*-algebras and
  second {\v{C}}ech cohomology}}, Indiana University Mathematics Journal
  \textbf{29} (1980), no.~6, 799--822.

\bibitem{phillips_tracial_2012}
N.~Christopher Phillips, \emph{The tracial {R}okhlin property is generic},
  {arXiv}:1209.3859 (2012).

\bibitem{popa_deformation_2007}
Sorin Popa, \emph{Deformation and rigidity for group actions and von {N}eumann
  algebras}, International Congress of Mathematicians. Vol. I, European
  Mathematical Society, Zürich, 2007, pp.~445--477.

\bibitem{rordam_simple_2003}
Mikael R{\o}rdam, \emph{A simple {C}*-algebra with a finite and an infinite
  projection}, Acta Mathematica \textbf{191} (2003), no.~1, 109--142.

\bibitem{rosenberg_algebraic_1994}
Jonathan Rosenberg, \emph{Algebraic {$K$}-theory and its applications},
  Graduate Texts in Mathematics, vol. 147, Springer-Verlag, New York, 1994.

\bibitem{rosendal_generic_2009}
Christian Rosendal, \emph{The generic isometry and measure preserving
  homeomorphism are conjugate to their powers}, Fundamenta Mathematicae
  \textbf{205} (2009), no.~1, 1{\textendash}27.

\bibitem{sabok_completeness_2013}
Marcin Sabok, \emph{Completeness of the isomorphism problem for separable
  {C}*-algebras}, {arXiv}:1306.1049 (2013).

\bibitem{sato_rohlin_2010}
Yasuhiko Sato, \emph{The {R}ohlin property for automorphisms of the
  {J}iang{\textendash}{S}u algebra}, Journal of Functional Analysis
  \textbf{259} (2010), no.~2, 453--476.

\bibitem{schochet_algebraic_1994}
Claude Schochet, \emph{Algebraic topology and {C}*-algebras}, {C}*-algebras:
  1943-1993 (San Antonio, {TX}, 1993), Contemporary Mathematics, vol. 167,
  Amer. Math. Soc., Providence, {RI}, 1994, pp.~218--231.

\bibitem{segal_c*-algebras_1994}
Irving Segal, \emph{{C}*-algebras and quantization}, {C}*-algebras: 1943-1993
  (San Antonio, {TX}, 1993), Contemporary Mathematics, vol. 167, Amer. Math.
  Soc., Providence, {RI}, 1994, pp.~54--65.

\bibitem{silver_counting_1980}
Jack~H. Silver, \emph{Counting the number of equivalence classes of {B}orel and
  coanalytic equivalence relations}, Annals of Mathematical Logic \textbf{18}
  (1980), no.~1, 1--28.

\bibitem{toms_independence_2005}
Andrew~S. Toms, \emph{On the independence of {K}-theory and stable rank for
  simple {C}*-algebras}, Journal f{\"{u}}r die reine und angewandte Mathematik
  \textbf{578} (2005), 185--199.

\bibitem{toms_classification_2008}
\bysame, \emph{On the classification problem for nuclear {C}*-algebras}, Annals
  of Mathematics \textbf{167} (2008), no.~3, 1029--1044.

\bibitem{toms_comparison_2009}
\bysame, \emph{Comparison theory and smooth minimal {C}*-dynamics},
  Communications in Mathematical Physics \textbf{289} (2009), no.~2, 401--433.

\bibitem{toms_strongly_2007}
Andrew~S. Toms and Wilhelm Winter, \emph{Strongly self-absorbing
  {C}*-algebras}, Transactions of the American Mathematical Society
  \textbf{359} (2007), no.~8, 3999--4029.

\bibitem{wassermann_tensor_1976}
Simon Wassermann, \emph{On tensor products of certain group {C}*-algebras},
  Journal of Functional Analysis \textbf{23} (1976), no.~3, 239--254.

\bibitem{wiersma_c*-norms_2014}
Matthew Wiersma, \emph{C*-norms for tensor products of discrete group
  {C}*-algebras}, {arXiv}:1406.2654 (2014).

\bibitem{zapletal_forcing_2014}
Jind{\v{r}}ich Zapletal, \emph{Forcing {B}orel reducibility invariants}, in
  preparation.

\bibitem{zapletal_analytic_2013}
\bysame, \emph{Analytic equivalence relations and the forcing method}, Bulletin
  of Symbolic Logic \textbf{19} (2013), no.~4, 473--490.

\bibitem{zielinski_complexity_2014}
Joseph Zielinski, \emph{The complexity of the homeomorphism relation between
  compact metric spaces}, {arXiv}:1409.5523 (2014).

\end{thebibliography}
\end{document}